\documentclass[11pt]{article}

\usepackage[margin=0.90in]{geometry}
\usepackage{amsmath,amssymb,amsthm}
\usepackage{array}
\usepackage{microtype}
\usepackage[hidelinks]{hyperref}

\newtheorem{theorem}{Theorem}
\newcommand{\Aut}{\operatorname{Aut}}

\title{A Ten-Vertex Counterexample to a Conjecture on Unstable Graphs}
\author{Prateek R. Srivastava\\
\href{mailto:prs7786@g.rit.edu}{\texttt{prs7786@g.rit.edu}}}
\date{14 August 2026}

\begin{document}
\maketitle

\begin{abstract}
Mizzi conjectured that every nontrivially unstable graph contains cycles
$C_k$ and $C_{2k}$ for some odd $k$.  We give a connected,
nonbipartite, vertex-determining counterexample on ten vertices.  Its
instability is certified by an explicit nontrivial two-fold automorphism,
and its complete set of simple-cycle lengths is $\{5,5,6\}$.
\end{abstract}

\section{The conjecture and the graph}

All graphs here are finite and simple.  For a graph $G$, write $N(v)$ for
the open neighbourhood of $v$.  A pair of permutations $(\alpha,\beta)$
of $V(G)$ is a \emph{two-fold automorphism} if
\[
 uv\in E(G)\quad\Longleftrightarrow\quad
 \alpha(u)\beta(v)\in E(G).
\]
Equivalently,
\begin{equation}\label{eq:tf}
                 \beta(N(v))=N(\alpha(v))
                 \qquad(v\in V(G)).
\end{equation}
This terminology is standard in the literature on unstable graphs
\cite{LauriMizziScapellato}.
It is nontrivial when $\alpha\ne\beta$.  A graph with such an
automorphism is unstable; this also follows directly by applying $\alpha$
and $\beta$ on the two colour classes of its canonical double cover
\cite[Theorem~2.1]{Mizzi}.

Mizzi conjectured that every TF-cousin pair and every nontrivially
unstable graph contains $C_k$ and $C_{2k}$ for some odd $k$
\cite[Section~7]{Mizzi}.  The following graph refutes the assertion about
unstable graphs.

\begin{theorem}\label{thm:main}
Let $G$ have vertex set $\{1,\ldots,10\}$ and edge set
\begin{equation}\label{eq:edges}
\{12,\,1\,10,\,24,\,34,\,36,\,45,\,57,\,5\,10,
  6\,10,\,79,\,8\,10\}.
\end{equation}
Then $G$ is connected and nontrivially unstable, but there is no odd $k$
for which $G$ contains both $C_k$ and $C_{2k}$.
\end{theorem}

\section{Verification}

The open neighbourhoods of $G$ are
\[
\begin{array}{c|cccccccccc}
v&1&2&3&4&5&6&7&8&9&10\\ \hline
N(v)&
2,10&1,4&4,6&2,3,5&4,7,10&3,10&5,9&10&7&1,5,6,8
\end{array}
\]
They are pairwise distinct.  The graph is connected, and the cycle
$1,2,4,5,10,1$ shows that it is nonbipartite.  Thus it satisfies the
paper's exclusions of the two classes called trivially unstable
\cite[Section~1]{Mizzi}.

Put
\[
                    \alpha=(1\ 6),\qquad \beta=(2\ 3).
\]
Both permutations fix every vertex not displayed.  Direct substitution in
the neighbourhood table verifies equation~\eqref{eq:tf} in all ten cases.
Hence $(\alpha,\beta)$ is a two-fold automorphism.  For completeness, the
induced map on the canonical double cover is
\[
 (v,0)\longmapsto(\alpha(v),0),\qquad
 (v,1)\longmapsto(\beta(v),1).
\]
Equation~\eqref{eq:tf} says exactly that this map preserves adjacency.
It preserves the two cover classes but does not act by the same
permutation on them.  It therefore does not belong to
$\Aut(G)\times\mathbb Z_2$, proving directly that $G$ is unstable.

It remains to determine all cycle lengths.  Vertices $8$ and $9$ are
leaves, and vertex $7$ lies only on the pendant path $5,7,9$; none can
belong to a cycle.  Deleting $7,8,9$ leaves a theta graph consisting of
the three internally disjoint paths from $4$ to $10$,
\[
4,2,1,10,\qquad 4,3,6,10,\qquad 4,5,10,
\]
of lengths $3,3,2$.  Every simple cycle in a theta graph is the union of
two of its three paths.  Consequently the three simple cycles of $G$ have
lengths
\[
                         3+3=6,\qquad 3+2=5,\qquad 3+2=5.
\]
For $k=3$ the graph has no $C_3$; for $k=5$ it has no $C_{10}$; and for
odd $k\ge7$, $C_{2k}$ has more than ten vertices.  This proves
Theorem~\ref{thm:main}.

\section*{Reproducibility and disclosure}

A dependency-free exact verifier reconstructs the graph from
\eqref{eq:edges}, checks all ten neighbourhood identities, checks the
induced automorphism of the canonical double cover, and enumerates every
simple cycle.  The example, verifier, and drafting were developed with
assistance from OpenAI Codex.  The proof above is self-contained, and the
author takes responsibility for all claims.

\begin{sloppypar}

\end{sloppypar}


\begin{thebibliography}{9}
\raggedright

\bibitem{LauriMizziScapellato}
J.~Lauri, R.~Mizzi, and R.~Scapellato,
\emph{Unstable graphs: a fresh outlook via TF-automorphisms},
Ars Math. Contemp. \textbf{8} (2015), 115--131.

\bibitem{Mizzi}
R.~Mizzi,
\emph{Lifting and Folding: A Framework for Unstable Graphs and
TF-Cousins},
\href{https://arxiv.org/abs/2603.27559v2}{arXiv:2603.27559v2}
[math.CO], 2026.

\end{thebibliography}
\end{document}